\documentclass[11pt,a4paper]{article}%

\usepackage{amsthm,amssymb}

\usepackage{authblk}
\theoremstyle{definition}
\newtheorem{definition}{Definition}[section]

\usepackage[T1]{fontenc}
\usepackage{graphicx}

\usepackage[utf8]{inputenc}
\usepackage{lmodern} %
\usepackage[english]{babel}

\usepackage{amsmath,amscd,latexsym}
\usepackage{bm}
\usepackage{hyperref}

\usepackage{algorithm}
\usepackage{algpseudocode}
\algrenewcommand\algorithmicensure{\textbf{Input:}}

\usepackage{color}
\usepackage{graphicx}
\graphicspath{{figures/}}
\usepackage{overpic}
\usepackage{subcaption}
\usepackage{placeins} %
\usepackage{siunitx}
\definecolor{shapes_1}{RGB}{230, 97, 1}
\definecolor{shapes_2}{RGB}{253, 140, 0}
\definecolor{shapes_3}{RGB}{178, 171, 210}
\definecolor{shapes_4}{RGB}{94, 60, 153}
\definecolor{shapes_5}{RGB}{43, 131, 186}

\newlength\figurewidth
\newlength\figureheight
\newlength\svgwidth

\newcommand{\R}{\mathbb{R}}

\newcommand{\dx}{\,\mathrm{d} \bm{x}}

\newcommand{\shapes}{\bm{u}}
\newcommand{\compdomain}{D_{\bm{u}}}
\newcommand{\velocity}{\bm{v}}
\newcommand{\pressure}{p}
\newcommand{\state}{\bm{y}}

\newcommand{\deformation}{\bm{W}}
\DeclareMathOperator{\vol}{vol}
\DeclareMathOperator{\bary}{bary}

\begin{document}
\title{A product shape manifold approach for optimizing piecewise-smooth shapes\thanks{This work has been partly supported by the German Research Foundation (DFG) within the priority program SPP~1962 under contract number WE~6629/1-1 and 
		by the state of Hamburg (Germany) within the Landesforschungsförderung under project %
		SENSUS with project number LFF-GK11.}}
 \author[1]{Lidiya Pryymak}
\author[2]{Tim Suchan}
\author[3]{Kathrin Welker}
\affil[1]{\small TU Bergakademie Freiberg, \texttt{lidiya.pryymak@math.tu-freiberg.de}
}
\affil[2]{Helmut Schmidt University, \texttt{suchan@hsu-hh.de}}
\affil[3]{TU Bergakademie Freiberg, \texttt{kathrin.welker@math.tu-freiberg.de}}

\date{}

\providecommand{\keywords}[1]
{
  \small	
  \textbf{\textit{Key words:}} #1
}
\providecommand{\ams}[1]
{
  \small	
  \textbf{\textit{AMS subject classifications:}} #1
}
\maketitle     
\vspace*{-1cm}         %
\begin{abstract}
Spaces where each element describes a shape, so-called shape spaces, are of particular interest in shape optimization and its applications. 
Theory and algorithms in shape optimization are often based on techniques from differential geometry.
Challenges arise when an application demands a non-smooth shape, which is commonly-encountered as an optimal shape for fluid-mechanical problems. 
In order to avoid the restriction to infinitely-smooth shapes of a commonly-used shape space, we construct a space containing shapes in $\R^2$ that can be identified with a Riemannian product manifold but at the same time admits piecewise-smooth curves as elements. 
We combine the new product manifold with an approach for optimizing multiple non-intersecting shapes.
For the newly-defined %
shapes, adjustments are made in the known shape optimization definitions and algorithms to ensure their usability in applications.
Numerical results regarding a fluid-mechanical problem constrained by the Navier-Stokes equations, where the viscous energy dissipation is minimized, show its applicability.
\end{abstract}


\smallskip
\noindent
\emph{\small\textbf{AMS subject classifications:}}
{\small
49Q10, 53C15, 58D10, 35Q30, 65K05}
\section{Introduction}

 Shape optimization is commonly-applied in enginnering in order to optimize shapes w.r.t. to an objective functional that relies on the solution of a partial differential equation (PDE). The PDE is required to model the underlying physical phenomenon, e.g. elastic displacements due to loadings or fluid movement due to pressure differences. Different methods are available for the shape optmization, however we focus on gradient-based techniques on shape spaces.

An ideal shape space would enable the usage of classical optimization methods like gradient descent algorithms. Since this is usually not the case, it is desirable to define a shape $u$ to be an element of a Riemannian manifold.
 An important example of a smooth\footnote{Throughout this paper, the term smooth shall refer to infinite differentiability.} manifold allowing a Riemannian structure is
the shape space 
$$B_e:=B_e(S^1, \mathbb{R}^2) := \text{Emb}(S^1, \mathbb{R}^2) /\text{Diff}(S^1).$$
An element of $B_e$ is a smooth simple closed curve in $\mathbb{R}^2$. The space was briefly investigated in~\cite{Michor2006}. The existence of Riemannian metrics, geodesics or, more generally, the differential geometric structure of $B_e$ (cf., e.g.~\cite{Michor2007,Michor2006}) reveals many possibilities like the computation of the shape gradient in shape optimization (cf., e.g.~\cite{Schulz2014}).
However, since an element of $B_e$ is a smooth curve in $\mathbb{R}^2$, the shape space is in general not sufficient to carry out optimization algorithms on piecewise-smooth shapes, which are often encountered as an optimal shape for fluid-mechanical problems, see e.g.~\cite{Pironneau1973} for a prominent example. In particular, we are interested in shapes with kinks. Such piecewise-smooth shapes are generally not elements of a shape space that provides the desired geometrical properties for applications in shape optimization. Some effort has been put into constructing a shape space that contains non-smooth shapes, however so far only a diffeological space structure could be found, cf. e.g.~\cite{Welker2016,Welker2021}. 
A further issue for many applications in shape optimization \cite{Albuquerque2020,Cheney1999,Kwon2002}, such as the electrical impedance tomography, is to consider multi-shapes. A first approach for optimizing  smooth multi-shapes has been presented in~\cite{Geiersbach2022}. 

In this paper, we aim to construct a novel shape space holding a Riemannian structure for optimizing  piecewise-smooth multi-shapes. %
The structure of the paper is as follows: In section \ref{sec:ShapeOpt}, we extend the findings related to multi-shapes
  in~\cite{Geiersbach2022} to a novel shape space considering piecewise-smooth shapes. Hereby, we use the fact that the space of simple, open curves $$B_e([0,1], \mathbb{R}^2) :=\text{Emb}([0,1], \mathbb{R}^2) /\text{Diff}([0,1])$$
is a smooth manifold as well (cf. \cite{Michor1980}) and interpret a closed curve with kinks as a glued-together curve of smooth, open curves, i.e.,  elements of $B_e([0,1], \mathbb{R}^2)$. Moreover, we derive a shape optimization procedure on the novel shape space.
In section \ref{sec:Numerics}, we apply the presented optimization technique to a shape optimization problem constrained by Navier-Stokes equations and present numerical results.

\label{sec:Introduction}

\section{Product space for optimizing piecewise-smooth shapes}
\label{sec:ShapeOpt}

In this section, we aim to construct a gradient descent algorithm for optimizing piecewise-smooth multi-shapes, e.g. the multi-shape $u=(u_1,u_2)$ from figure~\ref{fig:computational_domain_sketch}. In section \ref{sec:product_shape_space}, we therefore introduce a novel shape space which has the structure of a Riemannian product manifold. An optimization algorithm on the novel shape space is formulated in section \ref{sec:optimization_technique_on_product_manifold_of_piecewise_smooth_shapes}.

\subsection{Product shape space} \label{sec:product_shape_space}
In the following, we introduce a novel shape manifold, whose structure will later be used to optimize piecewise-smooth shapes.
The construction of the novel shape space is based on a Riemannian product manifold. Therefore, we first investigate the structure of product manifolds.

We define $(\mathcal{U}_i, G^i)$ to be Riemannian manifolds equipped with the Riemannian metrics $G^i$ for all $i=1,\ldots,N\in\mathbb{N}$.
The Riemannian metric $G^i$ at the point $p \in \mathcal{U}_i$ will be denoted by $$G_{p}^i(\cdot,\cdot)\colon T_{p} \mathcal{U}_i \times T_{p} \mathcal{U}_i \to \R,$$ 
where $T_{p} \mathcal{U}_i$ denotes the tangent space at a point $p \in \mathcal{U}_i$. We then define the product 
manifold as
$$  \mathcal{U}^N := \mathcal{U}_1 \times \ldots \times \mathcal{U}_N= \prod_{i=1}^{N}\mathcal{U}_i.$$
As shown in~\cite{Geiersbach2022}, for the tangent space of product manifolds it holds $$T_{\tilde{u}} \mathcal{U}^N \cong T_{\tilde{u}_1} \mathcal{U}_1 \times \dots \times  T_{\tilde{u}_N} \mathcal{U}_N.$$
Moreover, a product metric can be defined as
\begin{align}
\mathcal{G}^N = \sum\limits_{i=1}^N \pi_i^{*} G^i,
\label{eq:product_metric}
\end{align}
where $\pi_i^{*}$ are the pushforwards associated with canonical projections. 
It is obvious to use the space $B_e$ defined in section \ref{sec:Introduction} to construct a specific product shape space. An issue arises for non-smooth shapes, e.g. the shape $u_1$ from figure \ref{fig:computational_domain_sketch}.
To fix this issue, we now introduce the new multi-shape space for $s$ shapes built on the Riemannian product manifold $\mathcal{U}^N$.
\begin{definition} 
	\label{def:product_metric}
 Let $(\mathcal{U}_i, G^i)$ be Riemannian manifolds equipped with Riemannian metrics $G^i$ for all $i=1,\dots,N$. Moreover, $  \mathcal{U}^N :=\prod_{i=1}^{N}\mathcal{U}_i$. For $s \in \mathbb{N}$, we define the $s$-dimensional shape space on $\mathcal{U}^N$ by
\begin{align*}
M_s(\mathcal{U}^N):= \{ u=(u_1, \dots ,u_s) \mid \,& u_j \in  \prod\limits_{l=k_j}^{k_j+n_j-1} \mathcal{U}_l ,\, \sum\limits_{j=1}^s n_j = N \text{ and }\\
&k_1=1, \, k_{j+1}=k_j+n_j
\,\forall j =1,\dots,s-1  \}.
\end{align*}	
\end{definition}
With definition \ref{def:product_metric}, an element in $M_s(\mathcal{U}^N)$ is defined as a group of $s$ shapes $u_1,\dots,u_s$, where each shape $u_j$ is an element of the product of $n_j$ smooth manifolds.
 For $\mathcal{U}_l=B_e([0,1], \mathbb{R}^2)$ for $l=1, \dots, 12$ and $\mathcal{U}_{13}=B_e(S^1, \mathbb{R}^2)$, we can define the shapes presented in figure \ref{fig:computational_domain_sketch} by $(u_1,u_2)\in M_{2}(\mathcal{U}^{13})$, where $u_1 \in \prod_{l=1}^{12} \mathcal{U}_l$ and $u_2 \in \mathcal{U}_{13}$. 

For applications of definition \ref{def:product_metric} in shape optimization problems, it is of great interest to look at the tangent space of $M_s(\mathcal{U}^N)$. Since any element $u=(u_1,\dots,u_s)\in M_s(\mathcal{U}^N)$ can be understood as an element $\tilde{u}=(\tilde{u}_1,\dots, \tilde{u}_N)\in\mathcal{U}^N$, we set
$T_u M_s(\mathcal{U}^N) = T_{\tilde{u}} \mathcal{U}^N$ and 
$$G_u(\varphi,\psi)= G_{\tilde{u}}(\varphi, \psi) \ \forall \varphi, \psi \in T_u M_s(\mathcal{U}^N) = T_{\tilde{u}} \mathcal{U}^N.$$
Next, we consider shape optimization problems, i.e., we investigate so-called 
shape functionals.
A shape functional on $M_s(\mathcal{U}^N)$ is given by
$j\colon M_s(\mathcal{U}^N) \to \mathbb{R}$, $u\mapsto j(u)$.
In the following paragraph, we investigate solution techniques for shape optimization problems, i.e., for problems of the form
\begin{equation}
	\label{minproblem}
	\min_{u\in M_s(\mathcal{U}^N)} j(u).
\end{equation}

\subsection{Optimization technique on $M_s(\mathcal{U}^N)$ for optimizing piece\-wise-smooth shapes} \label{sec:optimization_technique_on_product_manifold_of_piecewise_smooth_shapes}

A theoretical framework for shape optimization depending on multi-shapes is presented in~\cite{Geiersbach2022}, where the optimization variable can be represented as a multi-shape belonging to a product shape space. 
Among other things, a multi-push\-forward and multi-shape gradient are defined; however, each shape is assumed to be an element of one shape space. In contrast, definition \ref{def:product_metric} also allows that a shape itself is represented by a product shape space. Therefore, we need to adapt the findings in~\cite{Geiersbach2022} to our setting.

To derive a gradient descent algorithm for a shape optimization problem as in \eqref{minproblem}, we need a definition for differentiating  a shape functional mapping from a smooth manifold to $\R$. For smooth manifolds, this is achieved using a pushforward.
\begin{definition} %
	\label{def:multi-pushforward}
For each shape $u \in M_s(\mathcal{U}^N)$, the multi-pushforward of a shape functional $j \colon M_s(\mathcal{U}^N) \rightarrow \mathbb{R}$ is given by the map
$$
(j_{*})_u \colon T_u M_s(\mathcal{U}^N) %
 \to \mathbb{R}, \ \varphi \mapsto \frac{\mathrm{d}}{\mathrm{d} t} j(\varphi(t))_{t=0}=(j \circ \varphi)^{'}(0).
$$
\end{definition}
Thanks to the multi-pushforward, we can  define the so-called multi-shape gradient, which is required for  optimization algorithms.
\begin{definition} %
	\label{def:multi-shapegradient}
The multi-shape gradient for a shape functional $j \colon M_s(\mathcal{U}^N) \rightarrow \mathbb{R}$ at the point $u \in M_s(\mathcal{U}^N)$ is given by $\psi \in T_uM_s(\mathcal{U}^N)$ satisfying	
$$
\mathcal{G}^N_{u}(\psi, \varphi)= (j_{*})_u \varphi \ \ \forall \varphi \in T_{\tilde{u}} M_s(\mathcal{U}^N).
$$
\end{definition}

We are now able to formulate a gradient descent algorithm on $M_s(\mathcal{U}^N)$ similar to the one presented in~\cite{Geiersbach2022}. %
For updating the multi-shape $u$ in each iteration, the multi-exponential map 
$$
\exp_{u}^N \colon T_u M_s(\mathcal{U}^N) \rightarrow M_s(\mathcal{U}^N), \ \varphi=(\varphi_1,\dots,\varphi_N) \mapsto (\exp_{\tilde{u}_1} \varphi_1,\dots, \exp_{\tilde{u}_N} \varphi_N)
$$
is used. The algorithm is depicted in algorithm~\ref{alg}.

\begin{algorithm}
	\caption{Gradient descent algorithm on $M_s(\mathcal{U}^N)$ with Armijo backtracking line search to solve \eqref{minproblem}}
	\label{alg}
	\begin{algorithmic}[1]
		\Ensure Initial shape $u=(u_1,\dots,u_s)=(\tilde{u}_1,\dots,\tilde{u}_N)= \tilde{u}$, 
		constants 
for Armijo backtracking and $\epsilon>0$ for break condition
		\While{$\| v \|_{\mathcal{G}^N} > \epsilon$}
		\State Compute the multi-shape gradient $v$ with respect to $\mathcal{G}^N$
		\State Compute stepsize $\alpha$ with Armijo backtracking 
		\State $u \leftarrow \exp^N_{u}(-\alpha v)$
		\EndWhile
	\end{algorithmic}
\end{algorithm}

So far, we have derived an optimization algorithm on $M_s(\mathcal{U}^N)$, i.e., an algorithm for optimizing a non-intersecting group of shapes, where each shape is an element of a product manifold with a varying number of factor spaces. With the main goal of this section in mind, we need to further restrict the choice of shapes in $M_s(\mathcal{U}^N)$ to glued-together piecewise-smooth shapes: We assume that $\mathcal{U}_i$ is either $B_e(S^1, \mathbb{R}^2)$ or $B_e([0,1], \mathbb{R}^2)$. 
Moreover, we assume that each shape $(u_1,\dots,u_s)$ is closed, where $u=(u_1,\dots,u_s)$ is chosen from $M_s(\mathcal{U}^N)$. By that we mean that if a shape is $u_j \in \prod\limits_{l=k_j}^{k_j+n_j-1} \mathcal{U}_l$, then either 
$$n_j=1 \ \text{and}  \ \mathcal{U}_{k_j}= B_e(S^1, \mathbb{R}^2)$$
or 
\begin{align*}
&\mathcal{U}_l= B_e([0,1], \mathbb{R}^2) \ \  \forall l=k_j,\dots,k_j+n_j-1 \text{ and for }
\\
 &u_j=(u_{k_j},\dots,u_{k_j+n_j-1})\text{, it holds that } 
\\
&u_{k_j+h}(1)=u_{k_j+h+1}(0) \ \ \forall h = 0,\dots,n_j-2 \ \ \text{and } u_{k_j}(0)=u_{k_j+n_j-1}(1).
\end{align*}

Finally, we want to address another important issue in shape optimization algorithms: the development of kinks in smooth shapes over the course of the optimization. If we view a smooth initial shape, e.g. $u_2$ from figure \ref{fig:computational_domain_sketch}, as an element in $B_e(S^1, \mathbb{R}^2)$ no kinks can arise during the optimization of the shape. 
An approach to fix this issue for applications, where the developments of kinks in shapes is desired, is to approximate a smooth shape with elements of $B_e([0,1], \mathbb{R}^2)$. A simple but sufficient choice is using initially straight lines connecting locations of possible kinks. In this manner, the multi-shape $u=(u_1,u_2)$ from figure \ref{fig:computational_domain_sketch} would be an element of  
\begin{equation}
\label{ShapeSpaceNum}
\begin{split}
&M_2(\mathcal{U}^{12+l_1+l_2})\text{, where } l_1,l_2\in\mathbb{N} \text{ and }\\
&u_1  \in \prod\limits_{l=1}^{12+l_1} \mathcal{U}_l= B_e([0,1], \mathbb{R}^2)^{12+l_1},\,
 u_2  \in \prod\limits_{l=13+l_1}^{12+l_1+l_2} \mathcal{U}_l= B_e([0,1], \mathbb{R}^2)^{l_2}.
\end{split}
\end{equation}

\section{Application to Navier-Stokes flow}
\label{sec:Numerics}

In the following, we apply algorithm~\ref{alg} to a shape optimization problem constrained by steady-state Navier-Stokes equations and geometrical constraints. %
In section \ref{sec:shape_derivative_etc}, we briefly describe the numerical implementation of algorithm~\ref{alg}. 
Afterwards, we formulate the optimization problem that will be considered for the numerical studies in section \ref{sec:model_formulation}, and finally, in section \ref{sec:numerical_results}, we describe the numerical results.

\subsection{Adjustments of algorithm \ref{alg} for numerical computations}\label{sec:shape_derivative_etc}
In order to ensure the numerical applicability of algorithm \ref{alg}, adjustments must be made. %
We define the space $\mathcal{W}:=\{ \deformation \in H^1(\compdomain, \R^2) \vert\, \deformation = \bm{0} \text{ on } \partial \compdomain \setminus \shapes \}$, and similarly to~\cite{Geiersbach2022}, we use an optimization approach based on partial shape derivatives, together with  the Steklov-Poincaré metric in equation~\eqref{eq:product_metric}. The Steklov-Poincaré metric is defined in \cite{SchulzSiebenbornWelker} and yields $\mathcal{G}^i(\bm{V}\vert_u,\deformation\vert_u) = a(\bm{V},\deformation)$ with a symmetric and coercive bilinear form $a\colon \mathcal{W} \times \mathcal{W}$. We replace the multi-shape gradient with the mesh deformation $\bm{V} \in \mathcal{W}$, which is obtained by replacing the multi-pushforward with the multi-shape derivative\footnote{We refer to \cite{Geiersbach2022} for the definition and details about the multi-shape derivative.} in definition \ref{def:multi-shapegradient}.
A common choice for the bilinear form when using the Steklov-Poincaré metric is linear elasticity
\begin{align}
	\int_{\compdomain} \bm{\varepsilon}(\bm{V}) : \bm{C} : \bm{\varepsilon}(\deformation) \dx = \mathrm{d} j(\bm{u}) [\deformation] \quad \forall \deformation \in \mathcal{W},
\end{align}
where $\bm{\varepsilon}(\bm{V})=\operatorname{sym}\operatorname{grad}(\bm{V})$ and $\bm{C}$ describes the linear elasticity tensor,  \mbox{$\bm{A} : \bm{B}$} is the standard Frobenius inner product and $ \mathrm{d} j(\bm{u}) [\deformation]$ denotes the shape derivative of $j$ at $\bm{u}$ in direction $\deformation$. 
Due to the equivalence of the Steklov-Poincaré metric and the bilinear form $a$, we replace the $\mathcal{G}^N$-norm in the stopping criterion of the algorithm~\ref{alg} with the $H^1$-norm in $D_{\shapes}$.
Finally, since the exponential map used in algorithm \ref{alg} is an expensive operation, it is common to replace it by a so-called retraction. In our computations, we use the retraction introduced in \cite{Schulz2018}. %

\subsection{Model formulation}  \label{sec:model_formulation}
\begin{figure}[tb]
	\centering
	\setlength\svgwidth{5cm}
\begingroup%
  \makeatletter%
  \providecommand\color[2][]{%
    \errmessage{(Inkscape) Color is used for the text in Inkscape, but the package 'color.sty' is not loaded}%
    \renewcommand\color[2][]{}%
  }%
  \providecommand\transparent[1]{%
    \errmessage{(Inkscape) Transparency is used (non-zero) for the text in Inkscape, but the package 'transparent.sty' is not loaded}%
    \renewcommand\transparent[1]{}%
  }%
  \providecommand\rotatebox[2]{#2}%
  \newcommand*\fsize{\dimexpr\f@size pt\relax}%
  \newcommand*\lineheight[1]{\fontsize{\fsize}{#1\fsize}\selectfont}%
  \ifx\svgwidth\undefined%
    \setlength{\unitlength}{230.20275287bp}%
    \ifx\svgscale\undefined%
      \relax%
    \else%
      \setlength{\unitlength}{\unitlength * \real{\svgscale}}%
    \fi%
  \else%
    \setlength{\unitlength}{\svgwidth}%
  \fi%
  \global\let\svgwidth\undefined%
  \global\let\svgscale\undefined%
  \makeatother%
  \begin{picture}(1,1.00000001)%
    \lineheight{1}%
    \setlength\tabcolsep{0pt}%
    \put(0,0){\includegraphics[width=\unitlength,page=1]{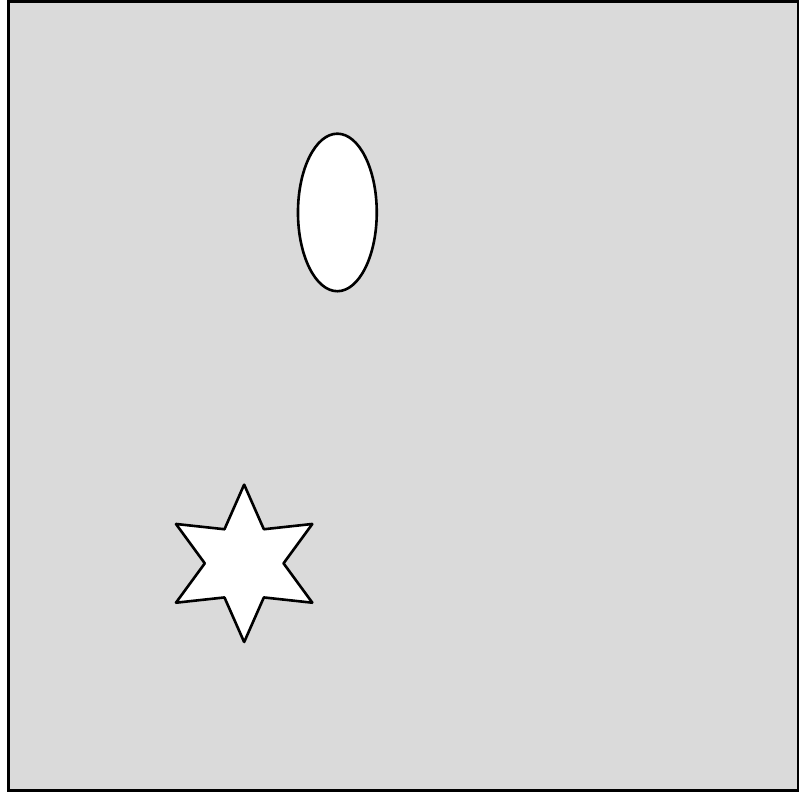}}%
    \put(0.24170821,0.37503696){\color[rgb]{0,0,0}\makebox(0,0)[t]{\lineheight{1.25}\smash{\begin{tabular}[t]{c}$u_1$\end{tabular}}}}%
    \put(0.32902976,0.80670375){\color[rgb]{0,0,0}\makebox(0,0)[t]{\lineheight{1.25}\smash{\begin{tabular}[t]{c}$u_2$\end{tabular}}}}%
    \put(0.71831187,0.45387213){\color[rgb]{0,0,0}\makebox(0,0)[t]{\lineheight{1.25}\smash{\begin{tabular}[t]{c}$\compdomain$\end{tabular}}}}%
    \put(0,0){\includegraphics[width=\unitlength,page=2]{computational_domain_sketch_svg-tex.pdf}}%
    \put(0.20448972,0.03779014){\color[rgb]{0,0,0}\makebox(0,0)[t]{\lineheight{1.25}\smash{\begin{tabular}[t]{c}$x_1$\end{tabular}}}}%
    \put(0.03079435,0.18792986){\color[rgb]{0,0,0}\makebox(0,0)[lt]{\lineheight{1.25}\smash{\begin{tabular}[t]{l}$x_2$\end{tabular}}}}%
  \end{picture}%
\endgroup%

	\caption{Sketch of two shapes $u_1$, $u_2$ surrounded by a  domain $\compdomain\subset \R^2$.}%
	\label{fig:computational_domain_sketch}
\end{figure}

We consider the problem
\begin{align}
	\label{eq:objectiveFunctional}
	\min_{\shapes \in  M_s(\mathcal{U}^N)}    j(\bm{u}) := \min_{\shapes \in  M_s(\mathcal{U}^N)} \int_{\compdomain} \frac{\mu}{2} \nabla \velocity : \nabla \velocity \dx,
\end{align}
where we constrain the optimization problem by the Navier-Stokes equations and choose $M_s(\mathcal{U}^N)$ as in (\ref{ShapeSpaceNum}). 
The state is denoted as $\state = (\velocity, \pressure)$ %
for which the Navier-Stokes equations can be found in standard literature and will be omitted here for brevity.
The material constants dynamic viscosity and density are defined as $\mu=1.81$ and $\rho=1.2\cdot 10^5$, respectively. We choose homogenous Dirichlet boundary conditions on the top and bottom boundary as well as on both shapes. The right boundary is modelled as homogenous Neumann, and the left boundary has the inhomogenous Dirichlet boundary condition $\velocity=(0.08421\, x_2 \, (x_2-1), 0)^\top$.
We choose the hold-all domain $D=(0,1)^2$, in which two shapes $u_1$ and $u_2$ are embedded as shown in figure~\ref{fig:computational_domain_sketch} with barycenters at $(0.3, 0.3)^\top$ and $(0.45, 0.75)^\top$, respectively. 

Additional geometrical constraints are required in order to avoid trivial solutions, see e.g.~\cite{Mohammadi2009,Mueller2021}, which are implemented as inequality constraints with an Augmented Lagrange approach as described in~\cite{Steck2018}. We restrict the area of each shape $\vol(u_i)$ to be at $100\%$ initial area. Further, the barycenter $\bary(u_1)$ is constrained to stay between $(-0.03, -0.05)^\top$ and $(0.04, 0.03)^\top$ of the initial position in $x$ and $y$ direction, respectively, and the barycenter $\bary(u_2)$ to stay between $(-0.075, -0.02)^\top$ and $(0.02, 0.05)^\top$ of the initial position.

\subsection{Numerical results}  \label{sec:numerical_results}

\begin{figure}[tb]
	\centering
	\setlength\figureheight{4.5cm} 
	\setlength\figurewidth{.45\textwidth}
	\includegraphics{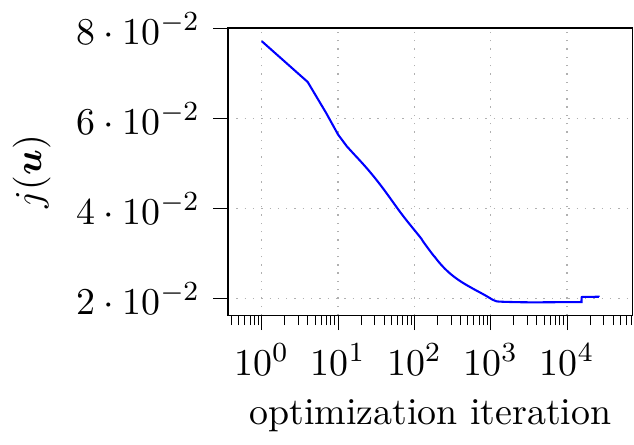}
	\hfill%
	\includegraphics{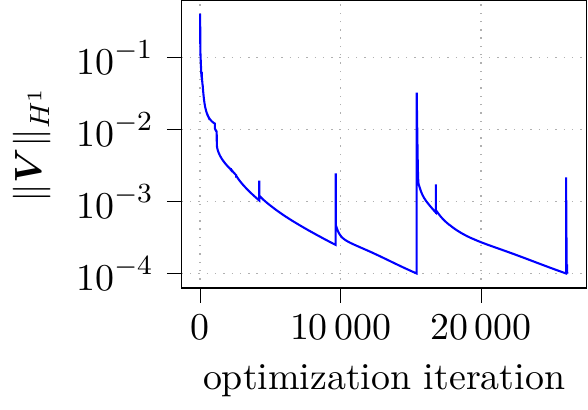}
	\caption{Optimization results: objective functional (left) and $H^1$-norm of the mesh deformation (right).}
\label{fig:NumericResults_objectiveFunctional_and_meshDeformNorm}
\end{figure}

\begin{figure}[htb]
	\centering
	{\captionsetup{justification=centering}
	\setlength\figureheight{.27\textwidth}  %
	\setlength\figurewidth{.27\textwidth}  %
	\hfill\includegraphics{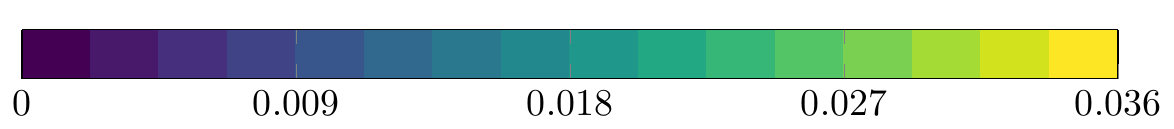}\\%
	\begin{subfigure}[t]{1.32\figurewidth}%
		\centering%
		\includegraphics{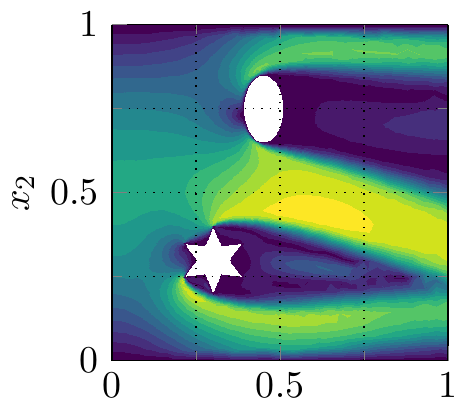}%
		\caption{Initial shapes:\\$j=0.07717$}%
		\label{fig:NumericResults_velocity_iter1}%
	\end{subfigure}\hspace*{.025\textwidth}%
	\begin{subfigure}[t]{1.08\figurewidth}%
		\centering%
		\includegraphics{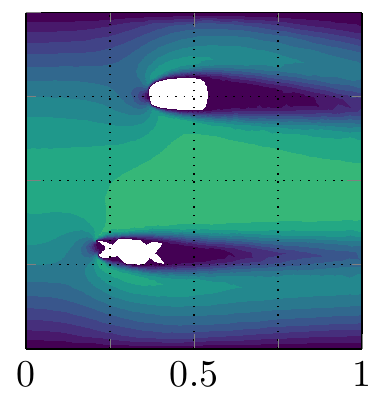}%
		\caption{Iteration 200:\\$j=0.02846$}%
		\label{fig:NumericResults_velocity_iter200}%
	\end{subfigure}\hspace*{.025\textwidth}%
	\begin{subfigure}[t]{1.08\figurewidth}%
		\centering%
		\includegraphics{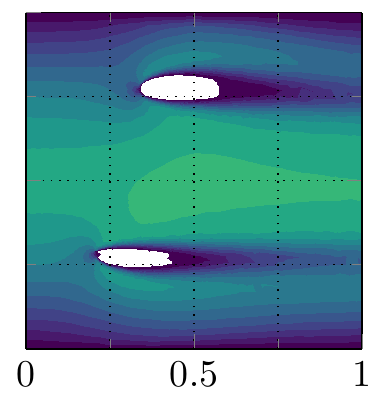}%
		\caption{Iteration 600:\\$j=0.02205$}%
		\label{fig:NumericResults_velocity_iter600}%
	\end{subfigure}\\%
	\begin{subfigure}[t]{1.32\figurewidth}%
		\centering%
		\includegraphics{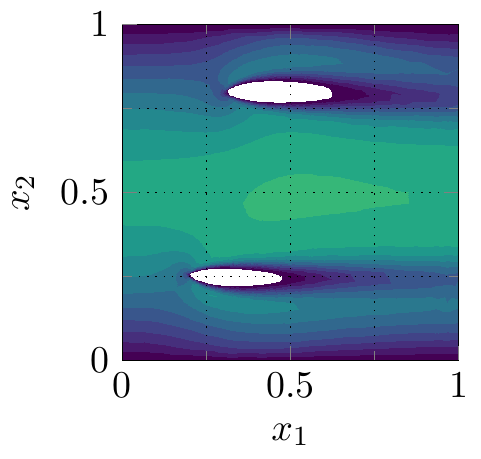}%
		\caption{Iteration 2\,000:\\$j=0.01918$}%
		\label{fig:NumericResults_velocity_iter2000}%
	\end{subfigure}\hspace*{.025\textwidth}%
	\begin{subfigure}[t]{1.08\figurewidth}%
		\centering%
		\includegraphics{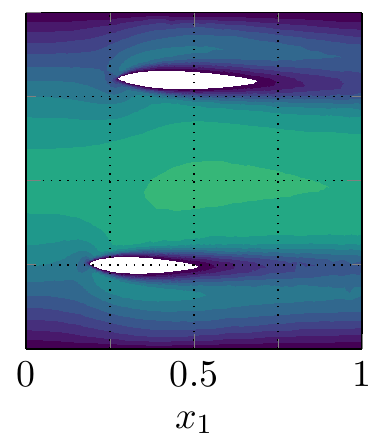}%
		\caption{Iteration 15\,000:\\$j=0.01917$}%
		\label{fig:NumericResults_velocity_iter15000}%
	\end{subfigure}\hspace*{.025\textwidth}%
	\begin{subfigure}[t]{1.08\figurewidth}%
		\centering%
		\includegraphics{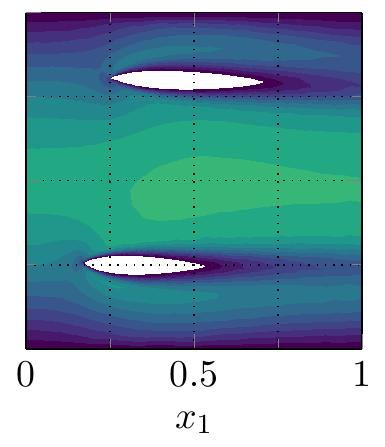}%
		\caption{Optimized shapes:\\$j=0.02034$}%
		\label{fig:NumericResults_velocity_iterEnd}%
	\end{subfigure}}%
	\vspace*{-1mm}
	\caption{Fluid velocity magnitude at different stages of the optimization. Figure~\ref{fig:NumericResults_velocity_iterEnd} has an increased objective functional value in comparison to figure~\ref{fig:NumericResults_velocity_iter2000} and~\ref{fig:NumericResults_velocity_iter15000}, however it fulfills the geometrical constraints while the others do not yet.}\vspace*{-3mm}
	\label{fig:NumericResults_velocity}
\end{figure}

The computational domain is discretized with 3512 nodes and 7026 triangular elements using Gmsh \cite{Geuzaine2009} with standard Taylor-Hood elements. An automatic remesher is available in case the mesh quality deteriorates below a threshold. The optimization is performed in FEniCS 2019.1.0 \cite{Alnaes2015}. We use a Newton solver and solve the linearized system of equations using MUMPS 5.5.1 \cite{Amestoy2001,Amestoy2006}.
Armijo backtracking is performed as described in algorithm~\ref{alg} with $\tilde{\alpha}=0.0125$, $\sigma=10^{-4}$ and $\tilde{\rho}=\frac{1}{10}$. The stopping criterion for each gradient descent is reached when the $H^1$-norm of the mesh deformation is at or below $10^{-4}$.
The objective functional and the $H^1$-norm of the mesh deformation over the course of the optimization are shown in figure~\ref{fig:NumericResults_objectiveFunctional_and_meshDeformNorm} and the magnitude of the fluid velocity in the computational domain before, during, and after optimization can be found in figure~\ref{fig:NumericResults_velocity}. The optimized shapes can be seen in figure~\ref{fig:NumericResults_velocity} on the right. Over the course of the optimization we observe a reduction of the objective functional by approximately $74\%$. The norm of the mesh deformation shows an exponential decrease, similar to a classical gradient descent algorithm. The peaks are caused by remeshing or by the adjustment of Augmented Lagrange parameters. Initially, the optimizer is mainly concerned with obtaining a approximate optimized shape, see figure~\ref{fig:NumericResults_velocity_iter200}--\ref{fig:NumericResults_velocity_iter2000}, while the exact fulfillment of geometrical constraints is less relevant. The later stages optimize small features like the leading and trailing edge of the shape, see figure~\ref{fig:NumericResults_velocity_iter15000}, any suboptimal kinks that were still remaining are removed, and in figure~\ref{fig:NumericResults_velocity_iterEnd} the geometrical constraints are fulfilled with an infeasibility of below $10^{-6}$ after $k=7$ Augmented Lagrange iterations.

\section{Conclusion}

A novel shape space $M_s(\mathcal{U}^N)$ that provides both, a Riemannian structure and a possibility to consider glued-together shapes (in particular, shapes with kinks) is introduced. Additionally, an optimization algorithm, based on findings from~\cite{Geiersbach2022}, is formulated. 
The new algorithm is applied to solve an optimization problem constrained by the Navier-Stokes equations with additional geometrical inequality constraints, where we have observed a strong reduction of the objective functional and convergence of the gradient descent on $M_s(\mathcal{U}^N)$ similar to a classical gradient descent algorithm.
Forthcoming research should involve an investigation of the development of the shapes' overlaps (glued-together points) over the course of the optimization. Moreover, convergence statements need to be investigated.

 \bibliographystyle{plainurl}
 \bibliography{literature}

\end{document}